\documentclass[12pt]{article}
\usepackage[all]{xy}
\usepackage{amsmath}
\usepackage{amsfonts}
\usepackage{amssymb}
\usepackage{amscd}
\usepackage{amsthm}
\usepackage{latexsym}
\usepackage{amsbsy}

\newtheorem{prop}{Proposition}[section]
\newtheorem{theorem}{Theorem}[section]
\newtheorem{corollary}{Corollary}[section]

\begin{document}
\date{}
\title
{ Cyclic Homology of Hopf Comodule Algebras and Hopf Module Coalgebras }
\author{R. Akbarpour and  M. Khalkhali \\Department of Mathematics 
, University of Western Ontario\\London ON, Canada N6A 5B7 \\
\texttt{akbarpur@uwo.ca \; masoud@uwo.ca}}
\maketitle

\begin{abstract}
In this paper we construct a cylindrical module $A \natural \mathcal{H}$ for an $\mathcal{H}$-comodule algebra
$A$, where the antipode of the Hopf algebra $\mathcal{H}$ is bijective. We show that the cyclic module associated to the diagonal
of $A \natural \mathcal{H}$ is isomorphic with the cyclic module of the crossed product algebra $A \rtimes \mathcal{H}$. This enables
us to derive a spectral sequence for the cyclic homology of the crossed product algebra. We also construct a cocylindrical module for Hopf
module coalgebras and establish a similar spectral sequence to compute the cyclic cohomology of crossed product coalgebras. 
\end{abstract}



\section{Introduction}

Getzler and Jones in~\cite{gj93} introduced a method to compute the cyclic homology of a crossed product algebra $ A \rtimes G $,
where $G$ is a group that acts on an algebra $A$ by automorphisms. This method is based on constructing a cylindrical module, $A \natural G$, and
showing that $ \Delta( A \natural G) \cong \mathsf{C}_{\bullet}
(A \rtimes G)$, where $\Delta$ is the diagonal and $\mathsf{C}_{\bullet}$ the cyclic module functor. Then by using 
the Eilenberg-Zilber theorem for cylindrical modules, they obtained a quasi-isomorphism of mixed 
complexes $\Delta(A \natural G) \cong Tot_{\bullet}{(A \natural G)}$, and a 
spectral sequence converging to $HC_{\bullet}(A \rtimes G)$. This spectral sequence was first obtained, by a different method, 
by Feigin and Tsygan~\cite{fs86}.
We used the the same method in~\cite{rm01,rm02} to generalize their work to Hopf module algebras and Hopf 
comodule coalgebras where the antipode of $\mathcal{H}$ is assumed to be bijective. In this paper we treat the two remaining cases, 
i.e., we derive  spectral sequences that converge to the cyclic (co)homology of the crossed product 
(co)algebra associated to a Hopf comodule
algebra or to a Hopf module coalgebra. The formulas in~\cite{rm01},~\cite{rm02} and the present paper cannot be deduced from 
each other and in each case it required a lot of effort to discover them.

In this paper we work over a fixed commutative ground ring $k$. By algebra we mean a unital associative $k$-algebra and 
algebra homomorphisms are unit preserving. Similar conventions apply to coalgebras and Hopf algebras. The unadorned tensor
product $\otimes$ means tensor product over $k$.

We denote the comultiplication of a coalgebra by $\Delta$, its counit by $\epsilon$ and the antipode of a Hopf algebra
by $S$. We use Sweedler's notation and write $\Delta h = h^{(0)} \otimes h^{(1)},\; \Delta^2 h = (1 \otimes \Delta)\Delta  h=
h^{(0)} \otimes h^{(1)} \otimes h^{(2)}$, etc., where summation is understood. We assume the antipodes of our Hopf algebras are bijective.
While this assumption may not be crucial for our definition of crossed product algebras and coalgebras, 
it plays an important role in our construction of cylindrical modules in Sections 2 and 3. 

If $\mathcal{H}$ is a Hopf algebra, by an $\mathcal{H}$-module we mean a module over the underlying algebra of $\mathcal{H}$
and by an $\mathcal{H}$-comodule we mean a comodule over the underlying coalgebra of $\mathcal{H}$. We use Sweedler's notation
for comodules. Thus if $\rho: M \rightarrow \mathcal{H} \otimes M$ is the structure map of a left $\mathcal{H}$-comodule
$M$, we write $\rho(m)=m^{(\bar{1})} \otimes m^{(\bar{0})}$, where summation is understood.

We say that an algebra $A$ is a left $\mathcal{H}$-comodule algebra if it is a left $\mathcal{H}$-comodule and the coaction 
$\alpha:A \rightarrow \mathcal{H} \otimes A$ is an algebra map.
So, we have, for all $a,b$ in $A$, 
\begin{eqnarray*}
& &\alpha (ab) = a^{(\bar{1})} b^{(\bar{1})} \otimes a^{(\bar{0})} b^{(\bar{0})}  \\
& &\alpha(1)= 1 \otimes 1.
\end{eqnarray*}

Given a left $\mathcal{H}$-comodule algebra $A$, we define the \textit{crossed product algebra $A \rtimes \mathcal{H}$}, 
to be the $k$-module $A \otimes \mathcal{H}$ with the following product

\[
(a \otimes g)(b \otimes h) = ab^{(\bar{0})} \otimes S^{-1}(b^{(\bar{1})})gb^{(\bar{2})}h.
\]
It is easy to check that this is an associative product with $1 \otimes 1$ as its unit.

Let $C$ be a coalgebra and $\mathcal{H}$ a Hopf algebra. We say that $C$ is a left $\mathcal{{\mathcal{H}}}$-module coalgebra if 
 $C$ is a left $\mathcal{H}$-module and
$$ 
\Delta (h \cdot c) = h^{(0)} \cdot c^{(0)} \otimes h^{(1)} \cdot c^{(1)} , 
\;\;\;\; \epsilon( h \cdot c) = \epsilon (h) \epsilon(c).
$$
This says that the module map $\mu  : \mathcal{{\mathcal{H}}} \otimes C \rightarrow C $ is a coalgebra map, where $\mathcal{{\mathcal{H}}} \otimes C$ has the  
tensor product coalgebra structure.\\

If $C$ is an $\mathcal{H}$-module coalgebra, we define the \textit{crossed product coalgebra} $C \! > \! \blacktriangleleft \! \mathcal{H} $, to be the $k$-module
$C \otimes \mathcal{H}$ with the following coproduct
\[
\Delta(a \otimes g) = (a^{(0)} \otimes g^{(1)}) \otimes (S^{-1}(g^{(0)}S(g^{(2)})) \cdot a^{(1)} \otimes g^{(3)}).
\]
One can check that this is a coassociative coproduct with counit $\epsilon = \epsilon_C \otimes \epsilon_{\mathcal{H}}$.
In fact,
\begin{eqnarray*}
& &((\Delta \otimes 1)\circ \Delta)(a \otimes g)=(a^{(0)} \otimes g^{(2)}) \otimes\\
& &\hspace{1cm}( S^{-1}(g^{(1)}S(g^{(3)})) \cdot a^{(1)} \otimes g^{(4)})\otimes 
( S^{-1}(g^{(0)}S(g^{(5)})) \cdot a^{(2)} \otimes g^{(6)}),
\end{eqnarray*}
and,
\begin{eqnarray*}
& &((  1 \otimes \Delta)\circ \Delta)(a \otimes g)=(a^{(0)} \otimes g^{(2)}) 
\otimes( S^{-1}(g^{(1)}S(g^{(3)})) \cdot a^{(1)} \otimes g^{(6)})\otimes\\
& &\hspace{3cm} ( S^{-1}(g^{(5)}S(g^{(7)})) S^{-1}(g^{(0)}S(g^{(4)})) \cdot a^{(2)} \otimes g^{(8)})=\\
& &(a^{(0)} \otimes g^{(2)}) \otimes
( S^{-1}(g^{(1)}S(g^{(3)})) \cdot a^{(1)} \otimes g^{(4)})\otimes 
( S^{-1}(g^{(0)}S(g^{(5)})) \cdot a^{(2)} \otimes g^{(6)}).
\end{eqnarray*}

We assume the reader is familiar with notions of cyclic homology theory, in particular with (co)cyclic modules, para(co)cyclic
modules and (co)cylindrical modules. They can all be found in~\cite{rm01,rm02,gj93,ll92}.
We recall the definition of the cocyclic module $\mathsf{C}^{\bullet}(C)$ associated to a coalgebra $C$. We have $\mathsf{C}^n(C)=C^{\otimes(n+1)}, n \ge0,$
with coface, codegeneracy and cyclic maps defined by
\begin{eqnarray*}
& &\partial^i(a_0,a_1,\dots,a_n)=(a_0,\dots,a_i^{(0)},a_i^{(1)},\dots,a_n), \;\; 0 \le i \le n,\\
& &\partial^{n+1}(a_0,a_1,\dots,a_n)=(a_0^{(1)},\dots,\dots,a_n,a_0^{(0)}), \\
& &\sigma^i(a_0,a_1,\dots,a_n)=\epsilon(a_{i})(a_0,\dots,a_{i-1},a_{i+1},\dots,a_n), \;\; 1 \le i < n,\\
& &\tau(a_0,a_1,\dots,a_n)=(a_1,a_2,\dots,a_0). 
\end{eqnarray*}

\section{The Cylindrical Module $\mathbf {A} \mathbf{\natural^{ \textbf{\tiny{}}}} \mathbf{\mathcal{H}}.$}
In this section we introduce the cylindrical module \[A \natural^{\textbf{}} \mathcal{H}= 
\{ \mathcal{H}^{\otimes(p+1)}\otimes A^{\otimes (q+1)} \}_{p,q \ge 0} \] where 
$A$ is an $\mathcal{H}$-comodule algebra and $\mathcal{H}$ is a Hopf algebra with bijective antipode. We define the
operators $\tau^{p,q} ,\partial^{p,q} ,\sigma^{p,q}$ and $\bar{\tau}^{p,q} ,\bar{\partial}^{p,q} ,\bar{\sigma}^{p,q}$ 
as follows:
\begin{eqnarray} \label{eq:opcy1} \notag
& &\tau^{p,q}( g_0 , \dots , g_p \mid a_0 , \dots , a_q)  
= (a_q^{\bar{(1)}} \cdot (g_0, \dots , g_p) \mid
 a_q^{\bar{(0)}} , a_0, \dots ,a_{q-1}), \\ \notag
& &\partial^{p,q}_i(g_0,\dots,g_p \mid a_0 , \dots , a_q)\\  
& & \hspace{1cm} = (g_0,\dots,g_p \mid a_0,\dots ,a_{i} a_{i+1},\dots , a_q),  \;\;\; 0 \le i <q,  \\ \notag
& &\partial^{p,q}_q(g_0,\dots,g_p \mid a_0 , \dots , a_q)  
= (a_q^{\bar{(1)}} \cdot (g_0, \dots , g_p) \mid
 a_q^{\bar{(0)}} a_0, \dots ,a_{q-1}), \\ \notag  
& &\sigma^{p,q}_i(g_0,\dots,g_p \mid a_0 , \dots , a_q)  \\ \notag
& & \hspace{1cm} =(g_0,\dots,g_p \mid a_0,\dots , a_i , 1 , a_{i+1},\dots , a_q), \;\;\; 0 \le i \le q, 
\end{eqnarray}
where $h \cdot (g_0,\dots,g_p)=(h^{(0)} \cdot g_0, \dots,h^{(p)} \cdot g_p)$ and 
$h \cdot g = h^{(0)} g S(h^{(1)})$ is the conjugation action of $\mathcal{H}$ on itself.
\begin{eqnarray} \label{eq:opcy2} \notag
& &\bar{\tau}^{p,q}( g_0 , \dots , g_p \mid a_0 , \dots , a_q)
=(S^{-1}(a_0^{\bar{(1)}} \dots a_q^{\bar{(1)}}) \cdot g_p,g_0, \dots , g_{p-1} \mid  a_0^{\bar{(0)}}, \dots
, a_{q}^{\bar{(0)}}), \\ \notag
& &\bar{\partial}^{p,q}_i(g_0,\dots,g_p \mid a_0 , \dots , a_q)\\ \notag
& &\hspace{1cm} =(g_0,\dots,g_i g_{i+1},\dots,g_p \mid a_0,\dots , a_q), \;\;\; 0 \le i < p,  \\ 
& &\bar{\partial}^{p,q}_p (g_0,\dots,g_p \mid a_0 , \dots , a_q)\\ \notag
& &\hspace{1cm}=((S^{-1}(a_0^{\bar{(1)}} \dots a_q^{\bar{(1)}}) \cdot g_p )g_0, \dots , g_{p-1} \mid 
 a_0^{\bar{(0)}}, \dots
, a_{q}^{\bar{(0)}}), \\ \notag
& &\bar{\sigma}^{p,q}_i (g_0,\dots,g_p \mid a_0 , \dots , a_q)\\ \notag
& &\hspace{1cm}=(g_0,\dots,g_i,1,g_{i+1},\dots,g_p \mid a_0,\dots , a_q), \;\;\; 0 \le i \le p.
\end{eqnarray} 
\begin{theorem}
$A \natural^{\tiny{}} \mathcal{H}$ with the operators defined in (\ref{eq:opcy1}),(\ref{eq:opcy2}) is a cylindrical module.  
\end{theorem}
\begin{proof}
We check only the commutativity of the cyclic operators and the cylindrical condition.\\

\item[$\bullet \quad$] $\tau^{p,q} \bar{\tau}^{p,q} = \bar{\tau}^{p,q} \tau^{p,q}$
\begin{eqnarray*}
& &(\tau^{p,q} \bar{\tau}^{p,q})(g_0,\dots,g_p \mid a_0,\dots ,a_q)\\
&=&\tau^{p,q}( S^{-1}(a_0^{\bar{(1)}} \dots a_q^{\bar{(1)}}) \cdot g_p , g_0,\dots,g_{p-1} 
\mid a_0^{\bar{(0)}}, \dots,a_q^{\bar{(0)}})\\
&=&(a_q^{\bar{(1)}} \cdot (S^{-1}(a_0^{\bar{(1)}} \dots a_q^{\bar{(2)}}) \cdot g_p, g_0, 
\dots g_{p-1}) \mid a_q^{\bar{(0)}} ,a_0^{\bar{(0)}},  \dots a_{q-1}^{\bar{(0)}})\\
&=&(a_q^{\bar{(1)} (1)} \cdot (S^{-1}( a_0^{\bar{(1)}}   \dots a_{q-1}^{\bar{(1)}} a_q^{\bar{(1)} (0)}) \cdot g_p), 
a_q^{\bar{(1)} (2)} \cdot (g_0, 
\dots g_{p-1}) \mid a_q^{\bar{(0)}} ,a_0^{\bar{(0)}}, \dots a_{q-1}^{\bar{(0)}})\\
&=& (S^{-1}(a_{0}^{\bar{(1)}} \dots a_{q-1}^{\bar{(1)}}) \cdot g_p, 
a_q^{\bar{(1)}} \cdot (g_0, 
\dots g_{p-1}) \mid a_q^{\bar{(0)}} ,a_0^{\bar{(0)}}, \dots a_{q-1}^{\bar{(0)}}).
\end{eqnarray*}
\begin{eqnarray*}
& &(\bar{\tau}^{p,q} \tau^{p,q})(g_0,\dots,g_p \mid a_0,\dots ,a_q)\\
&=& \bar{\tau}^{p,q}(a_q^{\bar{(1)}} \cdot (g_0, \dots , g_p) \mid a_q^{\bar{(0)}}, a_0, \dots,a_{q-1})\\
&=& \bar{\tau}^{p,q}(a_q^{\bar{(1)}(0)} \cdot (g_0, \dots , g_{p-1}),a_q^{\bar{(1)}(1)} \cdot g_p \mid
 a_q^{\bar{(0)}}, a_0, \dots,a_{q-1})\\
&=&(S^{-1}(a_{q}^{\bar{(1)} } a_{0}^{\bar{(1)}} \dots a_{q-1}^{\bar{(1)}}  ) \cdot (a_q^{\bar{(2)}(1)} \cdot  g_{p})
,a_q^{\bar{(2)}(0)} \cdot (g_0, \dots, g_{p-1}) \mid
a_q^{\bar{(0)}} ,a_0^{\bar{(0)}}, \dots a_{q-1}^{\bar{(0)}})\\
&=&(S^{-1}(a_{q}^{\bar{(1)}(2) } a_{0}^{\bar{(1)}} \dots a_{q-1}^{\bar{(1)}}  ) \cdot (a_q^{\bar{(1)}(1)} \cdot  g_{p})
,a_q^{\bar{(1)}(0)} \cdot (g_0, \dots, g_{p-1})
 \mid
a_q^{\bar{(0)}} ,a_0^{\bar{(0)}}, \dots a_{q-1}^{\bar{(0)}})\\
&=& (S^{-1}(a_{0}^{\bar{(1)}} \dots a_{q-1}^{\bar{(1)}}) \cdot g_p, 
a_q^{\bar{(1)}} \cdot (g_0, 
\dots g_{p-1}) \mid a_q^{\bar{(0)}} ,a_0^{\bar{(0)}}, \dots a_{q-1}^{\bar{(0)}}).
\end{eqnarray*}

\item[$\bullet \quad$] $(\bar{\tau}^{p,q})^{p+1} \;
(\tau^{p,q})^{q+1}=(\tau^{p,q})^{q+1} \;
(\bar{\tau}^{p,q})^{p+1}=id_{p,q}$\\

$(\tau^{p,q})^{q+1} (\bar{\tau}^{p,q})^{p+1} (g_0,\dots,g_p \mid a_0, \dots , a_q) $
\begin{eqnarray*}
&=&(\tau^{p,q})^{q+1} ( S^{-1}(a_0^{(\overline{p+1})} \dots a_q^{(\overline{p+1})} ) \cdot g_0 , \dots ,
S^{-1}(a_0^{\bar{(1)}} \dots a_q^{\bar{(1)}} ) \cdot g_p \mid a_0^{\bar{(0)}} , \dots , a_q^{\bar{(0)}} )\\
&=&((a_0^{(\overline{1})} \dots a_q^{(\overline{1})}) \cdot (S^{-1}(a_0^{(\overline{2})} \dots a_q^{(\overline{2})} ) \cdot g_0 , \dots ,
S^{-1}(a_0^{(\overline{p+2})} \dots a_q^{(\overline{p+2})} ) \cdot g_p) \mid a_0^{\bar{(0)}} , \dots , a_q^{\bar{(0)}} )\\
&=&(g_0,\dots,g_p \mid a_0, \dots , a_q).
\end{eqnarray*}
\end{proof}
Let $\Delta(A \natural \mathcal{H})$ denote the diagonal of $A \natural \mathcal{H}.$ Since 
$A \natural \mathcal{H}$ is cylindrical, we have:
\begin{corollary}$\Delta (A \natural^{} \mathcal{H})$ is a cyclic module.
\end{corollary}
Let $\mathsf{C}_{\bullet}(A \rtimes \mathcal{H})$ denote the cyclic module of the crossed product algebra $A \rtimes \mathcal{H}$.
\begin{theorem} \label{iso1}
We have an isomorphism of cyclic modules 
\[\Delta (A \natural^{} \mathcal{H}) \cong \mathsf{C}_{\bullet}(A \rtimes \mathcal{H}).\]
\end{theorem}
\begin{proof}
One can check that $\phi$ and $\psi$ define  isomorphisms, inverse to each other, between 
$\Delta (A \natural^{} \mathcal{H})$ and $\mathsf{C}_{\bullet}(A \rtimes \mathcal{H})$ where,
\begin{eqnarray*}
& &\phi ( a_0 \otimes g_0,\dots ,a_n \otimes g_n )=\\
& &(S^{-1}(a_1^{\bar{(1)}} \dots a_n^{\bar{(1)}}) \cdot g_0 ,S^{-1}(a_2^{\bar{(2)}} \dots a_n^{\bar{(2)}}) \cdot
g_1, \dots ,S^{-1}(a_{n}^{(\overline{n})}) \cdot g_{n-1},g_n \mid a_0,a_1^{\bar{(0)}} , \dots , a_n^{\bar{(0)}}),  
\end{eqnarray*}
\begin{eqnarray*}
& &\psi (g_0, \dots , g_n \mid a_0 , \dots ,a_n)=\\
& &(a_0 \otimes ( a_1^{\bar{(1)}} \dots a_n^{\bar{(n)}}) \cdot g_0,a_1^{\bar{(0)}} \otimes ( a_2^{\bar{(1)}} \dots 
a_n^{(\overline{n-1})}) \cdot g_1,
\dots, a_{n-1}^{\bar{(0)}} \otimes ( a_n^{\bar{(1)}} ) \cdot g_{n-1}, a_n^{\bar{(0)}} \otimes g_n). 
\end{eqnarray*}
\end{proof}
We define an action of $\mathcal{H}$ on the first column of $A \natural^{} \mathcal{H}, A_{\mathcal{H}}^{\natural^{}} =\{\mathcal{H} \otimes A^{\otimes(n+1)} \}_{n \ge 0},$
by 
\[ h \cdot ( g \mid a_0, \dots ,a_n) = ((S^{-1}(a_0^{\bar{(1)}} 
\dots a_n^{\bar{(1)}}) \cdot h^{(0)})g S^{-1}(h^{(1)}) \mid a_0^{\bar{(0)}}, \dots ,a_n^{\bar{(0)}}), \]
where $\cdot$ is the conjugation action of $\mathcal{H}$ on itself. 
We define $C^{\mathcal{H}^{}}_{\bullet} (A)$ as the coinvariant space of $A^{\natural^{}}_{\mathcal{H}}$ under the above action,
so that 
\[ C^{\mathcal{H}^{}}_{\bullet} (A)=A_{\mathcal{H}}^{\natural^{}}/span\{h \cdot x - \epsilon(h)x |\; 
h \in \mathcal{H}, x \in A_{\mathcal{H}}^{\natural^{}}\}. \]

In order to compute the $\mathsf{E}^1$ and $\mathsf{E}^2$ term of the spectral sequence associated to $A \natural \mathcal{H}$, we
replace $A \natural \mathcal{H}$ by another cylindrical module.

Let $\mathsf{C}_q(C_{\mathcal{H}}^{\natural^{}})=\mathcal{H} \otimes C^{\otimes (q+1)}$ and
$\mathsf{C}_p(\mathcal{H},\mathsf{C}_q(C_{\mathcal{H}}^{\natural^{}}))= \mathcal{H}^{\otimes p} \otimes 
\mathsf{C}_q(C_{\mathcal{H}}^{\natural^{}}).$ We replace the cylindrical module  
$(A \natural^{} \mathcal{H},(\partial,\sigma,\tau),
(\bar{\partial},\bar{\sigma},\bar{\tau}))$ with a new cylindrical module 
$(\mathsf{C}(\mathcal{H},\mathsf{C}( A^{\natural^{}}_{\mathcal{H}})),
( \mathfrak{d},\mathfrak{s},\mathfrak{t}),(\bar{\mathfrak{d}},\bar{\mathfrak{s}},\bar{\mathfrak{t}}))$,
via isomorphisms defined by 
the maps $\beta : (A \natural^{} \mathcal{H})_{p,q} \rightarrow \mathsf{C}_p(\mathcal{H},\mathsf{C}_q(A_{\mathcal{H}}^{\natural^{}}))$
and $\gamma :  \mathsf{C}_p(\mathcal{H},\mathsf{C}_q(A_{\mathcal{H}}^{\natural^{}})) \rightarrow  (A \natural^{} \mathcal{H})_{p,q}$
\begin{eqnarray*}
& &\beta ( g_0, \dots , g_p \mid a_0 , \dots , a_q ) = ( g_1^{(0)}, \dots , g_p^{(0)} \mid g_0 g_1^{(1)} 
\dots g_p^{(1)} \mid
a_0, \dots , a_q ), \\
& &\gamma ( g_1 , \dots , g_p \mid g \mid a_0 , \dots , a_q) = (g S^{-1}(g_1^{(1)} \dots g_p^{(1)}),g_1^{(0)}, 
\dots , g_p^{(0)}  \mid 
a_0, \dots , a_q ).
\end{eqnarray*}

It can be shown that $\beta \circ \gamma = \gamma \circ \beta = id $. We find the operators,
$( \mathfrak{d},\mathfrak{s},\mathfrak{t}),(\bar{\mathfrak{d}},\bar{\mathfrak{s}},\bar{\mathfrak{t}})$
under this isomorphism. \\

$ \hspace{0.2cm} \bullet \quad $ 
First we compute $\bar{\mathfrak{b}} = \beta \bar{b} \gamma$. Since 
\begin{eqnarray*}
\bar{b} ( g_0, \dots , g_p \mid a_0 , \dots , a_q ) &=& \sum_{0 \le i < p} (-1)^i (g_0 , \dots , g_i g_{i+1} ,\dots ,
 g_p \mid a_0 , \dots , a_q ) \\
&+& (-1)^p ((S^{-1}(a_0^{\bar{(1)}} \dots a_q^{\bar{(1)}}) \cdot g_p )g_0, \dots , g_{p-1} \mid  a_0^{\bar{(0)}}, \dots
, a_{q}^{\bar{(0)}}), 
\end{eqnarray*}
we see that,
\begin{eqnarray*}
& &\bar{\mathfrak{b}}( g_1 , \dots , g_p \mid g \mid a_0 , \dots , a_q) = \beta \bar{b} (g S^{-1}(g_1^{(1)} \dots g_p^{(1)}),
g_1^{(0)}, \dots g_p^{(0)} \mid a_0 , \dots , a_q )\\
&=& ( g_2^{(0)}, \dots , g_p^{(0)} \mid \epsilon(g_1) g S^{-1}(g_2^{(2)} \dots g_p^{(2)}) g_2^{(1)} \dots  g_p^{(1)} \mid a_0 , \dots , a_q)\\ 
&+& \sum _{1 \le i < p} (-1)^i ( g_1^{(0)},\dots , g_i^{(0)}g_{i+1}^{(0)} , \dots , g_p^{(0)} \mid g S^{-1}(g_1^{(2)} \dots g_p^{(2)})g_1^{(1)} \dots g_p^{(1)} \mid a_0 , \dots , a_q)\\ 
&+& (-1)^p (g_1^{(0)}, \dots , g_{p-1}^{(0)} \mid (S^{-1}(a_0^{\bar{(1)}} \dots a_q^{\bar{(1)}}) 
\cdot g_p^{(0)} )g S^{-1}(g_1^{(2)} \dots g_{p-1}^{(2)} g_{p}^{(1)})g_1^{(1)} \dots g_{p-1}^{(1)}\\  
& &\hspace{10cm} \mid  a_0^{\bar{(0)}}, \dots
, a_{q}^{\bar{(0)}})
\end{eqnarray*}
\begin{eqnarray*}
&=& \epsilon(g_1)( g_2, \dots , g_p \mid  g  \mid a_0 , \dots , a_q)\\ 
&+& \sum _{0 \le i < p} (-1)^i ( g_1,\dots , g_i g_{i+1} , \dots , g_p  \mid g  \mid a_0 , \dots , a_q)\\ 
&+& (-1)^p (g_1, \dots , g_{p-1} \mid (S^{-1}(a_0^{\bar{(1)}} \dots a_q^{\bar{(1)}}) \cdot g_p^{(0)} )g S^{-1}(g_p^{(1)})  \mid  a_0^{\bar{(0)}}, \dots
, a_{q}^{\bar{(0)}})\\
&=& \delta ( g_1 , \dots , g_p \mid g \mid a_0 , \dots , a_q).
\end{eqnarray*}
Here the boundary $\delta$ is the Hopf-module boundary defined by
\begin{eqnarray*} 
& & \delta (g_1,g_2,\dots , g_p,m)= \epsilon (g_1) (g_2,\dots ,g_p,m) +\\
& &\sum_{i=1}^{p-1} (-1)^i (g_1 , \dots ,g_i g_{i+1},\dots , g_p ,m) + (-1)^p (g_1,\dots,g_{p-1},g_p \cdot m), 
\end{eqnarray*}
where $m \in M$ and $M$ is an $\mathcal{H}$-module.
 We also conclude that
\begin{eqnarray*}
& &\bar{\mathfrak{d}}_0 (g_1, \dots ,g_p \mid g \mid a_0, \dots , a_q) = 
\epsilon(g_1)(g_2 , \dots , g_p \mid g \mid a_0 , \dots , a_q),\\
& &\bar{\mathfrak{d}}_i (g_1, \dots ,g_p \mid g \mid a_0, \dots , a_q) = 
(g_1 , \dots , g_i g_{i+1} ,\dots , g_p \mid g \mid a_0 , \dots , a_q),\\
& &\bar{\mathfrak{d}}_p (g_1, \dots ,g_p \mid g \mid a_0, \dots , a_q) = 
(g_1 , \dots , g_{p-1} \mid g_p \cdot ( g \mid a_0 , \dots ,  a_q)).
\end{eqnarray*}
$\hspace{0.2cm} \bullet \quad$  For 
$\mathfrak{\bar{s}}_i = \beta \bar{\sigma}_i \gamma$ and $\mathfrak{\bar{t}} = \beta \bar{\tau} \gamma$ we have 
\begin{eqnarray*}
& &\mathfrak{\bar{s}}_i (g_1,\dots,g_p \mid g \mid a_0, \dots ,a_q) \\
& &=(g_1, \dots ,g_i,1,g_{i+1},\dots, g_p \mid g  \mid  a_0 , 
\dots, a_q)\;\; 0 \le i \le p,\\
& & \mathfrak{\bar{t}} (g_1, \dots , g_p \mid g \mid a_0, \dots ,a_q)
 \\  
& &= (g_1^{(0)}, \dots,g_p^{(0)} \mid 
(S^{-1}(a_0^{\bar{(1)}} \dots a_q^{\bar{(1)}}) \cdot (g S^{-1}(g_1^{(2)} \dots g_p^{(2)}))g_1^{(1)} \dots g_p^{(1)} \mid 
a_0^{\bar{(0)}} \dots a_q^{\bar{(0)}}) .
\end{eqnarray*}

$ \hspace{0.2cm} \bullet \quad $ Next we compute the operator  $\mathfrak{b} = \beta b \gamma$. We have 
\begin{eqnarray*}
& &\mathfrak{b}(g_1, \dots , g_p \mid g \mid a_0 , \dots , a_q) = \beta b (g S^{-1}(g_1^{(1)} \dots g_p^{(1)}),
g_1^{(0)}, \dots ,g_p^{(0)} \mid a_0 , \dots , a_q )\\
&=&  \sum_{ 0 \le i \le q}(-1)^i ( 
g_1, \dots ,g_p \mid g \mid a_0 , \dots ,a_{i} a_{i+1}, \dots , a_q ) \\
&+&(-1)^q (a_q^{(\overline{3p-1})} \cdot g_1^{(0)}, \dots,a_q^{\bar{(2)}} \cdot g_p^{(0)} \mid (a_q^{(\overline{3p+1})} 
\cdot (gS^{-1}(g_1^{(2)} \dots g_p^{(2)})))\\
& & (a_q^{(\overline{3p})}g_1^{(1)}S(a_q^{(\overline{3p-2})})a_q^{(\overline{3p-3})}g_2^{(1)} S(a_q^{(\overline{3p-5})})
\dots S(a_q^{\bar{(4)}}) a_q^{\bar{(3)}} g_p^{(1)}S(a_q^{\bar{(1)}}))\mid a_q^{\bar{(0)}} a_0 , \dots , a_{q-1})\\
&=&  \sum_{ 0 \le i \le q}(-1)^i ( 
g_1, \dots ,g_p \mid g \mid a_0 , \dots ,a_{i} a_{i+1}, \dots , a_q ) \\
&+&(-1)^q (a_q^{(\overline{2})} \cdot (g_1^{(0)}, \dots,g_p^{(0)}) \mid (a_q^{(\overline{4})} 
\cdot (gS^{-1}(g_1^{(2)} \dots g_p^{(2)})))\\
& &\hspace{7cm}(a_q^{(\overline{3})}g_1^{(1)} \dots g_p^{(1)}S(a_q^{\bar{(1)}}))\mid a_q^{\bar{(0)}} a_0 , \dots , a_{q-1})
\end{eqnarray*}
\begin{eqnarray*}
&=&  \sum_{ 0 \le i \le q}(-1)^i ( 
g_1, \dots ,g_p \mid g \mid a_0 , \dots ,a_{i} a_{i+1}, \dots , a_q ) \\
&+&(-1)^q (a_q^{(\overline{2})} \cdot (g_1, \dots,g_p) \mid a_q^{(\overline{1})} 
\cdot g \mid a_q^{\bar{(0)}} a_0 , \dots , a_{q-1}).
\end{eqnarray*}
So we conclude that for $ 0 \le i < q $,
\begin{eqnarray*}
& &\mathfrak{d}_i(g_1,\dots,g_p \mid g \mid a_0,\dots,a_q)=(g_1 , \dots , 
g_p \mid g \mid a_0, \dots,a_{i} a_{i+1},\dots,a_{q}),\\
& &\mathfrak{d}_q(g_1,\dots,g_p \mid g \mid a_0,\dots,a_q)= (a_q^{(\overline{2})} 
\cdot (g_1, \dots,g_p) \mid a_q^{(\overline{1})} 
\cdot g \mid a_q^{\bar{(0)}} a_0 , \dots , a_{q-1}).
\end{eqnarray*}
$\hspace{0.2cm} \bullet \quad$  We consider 
$\mathfrak{s}_i = \beta \sigma_i \gamma$ and $\mathfrak{t} = \beta \tau \gamma$. We have 
\begin{eqnarray*}
& &\mathfrak{s}_i (g_1,\dots,g_p \mid g \mid a_0, \dots ,a_q) \\
& &= \beta (g S^{-1} (g_1^{(1)} \dots g_p^{(1)}), g_1^{(0)}, \dots , g_p^{(0)} \mid  a_0 , \dots ,a_i,1, a_{i+1},\dots, a_q )\\ 
& &=(g_1, \dots , g_p \mid g  \mid  a_0 , 
\dots ,a_i,1, a_{i+1},\dots, a_q),\\
& & \mathfrak{t} (g_1, \dots , g_p \mid g \mid a_0. \dots , a_q ) \\
& &= \beta (a_q^{\bar{(1)}} \cdot (g S^{-1} (g_1^{(1)} \dots g_p^{(1)}), g_1^{(0)}, \dots , g_p^{(0)}) \mid a_q^{\bar{(0)}}, 
a_0 ,\dots, a_{q-1} )\\   
& &= (a_q^{(\overline{2})} \cdot (g_1, \dots,g_p) \mid a_q^{(\overline{1})} 
\cdot g \mid a_q^{\bar{(0)}}, a_0 , \dots , a_{q-1}).
\end{eqnarray*}

By the above computations we can state the following theorems.
\begin{theorem} \label{th:iso1}
The complex $(\mathsf{C} (A \natural^{} \mathcal{H}) \boxtimes \mathsf{W}, b + \mathbf{u} B , \bar{b} + \mathbf{u} \bar{B}) $ 
is isomorphic to the complex 
$(\mathsf{C}(\mathcal{H}, \mathsf{C}(A^{\natural}_{\mathcal{H}}) \boxtimes \mathsf{W}), 
\mathfrak{b} + \mathbf{u} \mathfrak{B} , \bar{\mathfrak{b}} + \mathbf{u} \bar{\mathfrak{B}}), $
where $\mathfrak{\bar{b}}$ is the Hopf-module boundary. 
($\mathsf{W}$ is defined in~\cite{gj93}). 
\end{theorem}
\begin{theorem} \label{th:iso2}
For  $p \ge 0$, $H_p(\mathcal{H},\mathsf{C}_{\bullet} (A^{\natural^{}}_{\mathcal{H}}))$, with the operators defined as follows,
are cyclic modules: 
\begin {eqnarray*} \label{eq:equ}
& &\mathfrak{t} (g_1, \dots ,g_p \mid g \mid a_0 , \dots , a_q )=(a_q^{(\overline{2})} \cdot 
(g_1, \dots,g_p) \mid a_q^{(\overline{1})} 
\cdot g \mid a_q^{\bar{(0)}}, a_0 , \dots , a_{q-1}),\\
& &\mathfrak{d}_i(g_1,\dots,g_p \mid g \mid a_0,\dots,a_q)=(g_1 , \dots , g_p \mid g \mid a_0,\dots, 
a_{i}a_{i+1} ,\dots,a_{q}), \;\; 0 \le i < q,\\
& &\mathfrak{d}_q(g_1,\dots,g_p \mid g \mid a_0,\dots,a_q)= 
(a_q^{(\overline{2})} \cdot (g_1, \dots,g_p) \mid a_q^{(\overline{1})} 
\cdot g \mid a_q^{\bar{(0)}} a_0 , \dots , a_{q-1}),\\ 
& &\mathfrak{s}_i (g_1, \dots ,g_p \mid g \mid a_0 , \dots , a_q )=(g_1, \dots , g_p \mid g \mid  a_0 , \dots ,
a_i,1, a_{i+1} ,\dots, a_q),\;\; 1 \le i \le q.
\end{eqnarray*}
\end{theorem}

\begin{proof} We check this only for $p=0.$ The general case is similar.
The operators are well defined on the coinvariant space. For example for $\mathfrak{t},$ since $ \mathfrak{t} \mathfrak{\bar{b}}= \mathfrak{\bar{b}} \mathfrak{t} $ and 
$$\mathfrak{\bar{b}}(h \mid g \mid a_0,\dots,a_n) = h \cdot (g \mid a_0,\dots,a_n) - \epsilon(h) (g \mid a_0,\dots,a_n),$$
we have
\begin{eqnarray*}
& &\mathfrak{t}(h \cdot (g \mid a_0,\dots,a_n) - \epsilon(h) (g \mid a_0,\dots,a_n))\\
& &=\mathfrak{t}\mathfrak{\bar{b}}(h \mid g \mid a_0,\dots,a_n)
= \mathfrak{\bar{b}}\mathfrak{t}(h \mid g \mid a_0,\dots,a_n).
\end{eqnarray*}
So $\mathfrak{t}$ is well defined on the coinvariant space.

To show that $\mathfrak{t}^{n+1} = 1,$ one can check  that in the coinvariant space
\begin{eqnarray*}
& &(g \mid a_0, \dots ,a_n) \equiv S(a_0^{\bar{(1)}} \dots a_n^{\bar{(1)}})\cdot(g \mid a_0^{\bar{(0)}}, \dots , a_n^{\bar{(0)}})\\
& &  \equiv (S(a_0^{\bar{(1)}} \dots a_n^{\bar{(1)}})\cdot g \mid a_0^{\bar{(0)}}, \dots , a_n^{\bar{(0)}}).
\end{eqnarray*}
Hence \[\mathfrak{t}^{n+1}(g \mid a_0, \dots ,a_n) =\mathfrak{t}^{n+1}( S(a_0^{\bar{(1)}} 
\dots a_n^{\bar{(1)}})\cdot g \mid a_0^{\bar{(0)}}, \dots ,a_n^{\bar{(0)}})\]
$\hspace{5cm} =(g \mid a_0, \dots ,a_n).$

\end{proof}
\begin{corollary} $C_{\bullet}^{\mathcal{H}^{\tiny{}}}(A)$ with the following operators is a cyclic module. 
\begin {eqnarray*} 
& &\mathfrak{t} (g \mid a_0 , \dots , a_n )=
(a_n^{\bar{(1)}} \cdot g  \mid a_n^{\bar{(0)}}, a_0 ,\dots, a_{n-1} ),\\   
& &\mathfrak{d}_i( g \mid a_0,\dots,a_n)=(g \mid a_0,\dots, a_{i}a_{i+1} ,\dots,a_{n}), \;\; 0 \le i < n,\\
& &\mathfrak{d}_n( g \mid a_0,\dots,a_n)= (a_n^{\bar{(1)}} \cdot g  \mid a_n^{\bar{(0)}} a_0 ,\dots, a_{n-1} ),\\   
& &\mathfrak{s}_i (g \mid a_0 , \dots , a_n )=( g \mid  a_0 , \dots ,a_i,1, a_{i+1} ,\dots, a_n), 0 \le i \le n.
\end{eqnarray*}
\end{corollary}

Now we use the Eilenberg-Zilber theorem for cylindrical modules, combined with Theorem~\ref{iso1} to conclude 
$$H_{\bullet}(Tot(A \natural^{} \mathcal{H}) \boxtimes \mathsf{W}) \cong HC_{\bullet} (\Delta( A \natural_{}  \mathcal{H});
\mathsf{W}) \cong HC_{\bullet} ( A \rtimes  \mathcal{H};
\mathsf{W}).$$
To compute the homology of the mixed complex 
$(Tot(A \natural^{} \mathcal{H}), b + \bar{b} + \mathbf{u} (B + \bar{B})),$
we filter it by the subspaces \[
\mathsf{F}_i Tot_{n}((A \natural^{} \mathcal{H} ) \boxtimes \mathsf{W}) = \sum_{q \le i,\;p+q=n} (\mathcal{H}^{\otimes(p+1)} \otimes A^{\otimes(q+1)}) \boxtimes \mathsf{W}. \]
If we separate the operator $b+ \bar{b} + \mathbf{u}(B + T \bar{B})$ as $\bar{b}+(b+  
\mathbf{u}B) + \mathbf{u}T \bar{B}$, from Theorems (\ref{th:iso1}) and (\ref{th:iso2}) we can deduce the 
following theorem for the spectral sequence associated to $Tot(A \natural \mathcal{H})$.

\begin{theorem}
There exists an spectral sequence that converges to $HC_{\bullet} (A \rtimes \mathcal{H}; \mathsf{W})$. The $\mathsf{E}^0$-term of this spectral sequence is isomorphic to the complex 
\[ \mathsf{E}^0_{pq}=(\mathsf{C}_q(\mathcal{H},
\mathsf{C}_p(A^{\natural^{}}_{\mathcal{H}}) \boxtimes \mathsf{W} ),\delta), \]
and the $\mathsf{E}^1$-term is  \[
\mathsf{E}^1_{pq}= (H_q( \mathcal{H} , \mathsf{C}_p(A^{\natural^{}}_{\mathcal{H}}) \boxtimes \mathsf{W} ) , 
\mathfrak{b} + \mathbf{u} \mathfrak{B}). \]
The $\mathsf{E}^2$-term of the spectral sequence is
\[ \mathsf{E}^2_{pq} = HC_p( H_q(\mathcal{H}, \mathsf{C}_{\bullet}(A^{\natural^{}}_{\mathcal{H}}));\mathsf{W}), \]
the cyclic homology of the cyclic module $H_q(\mathcal{H},\mathsf{C}_{\bullet}(A^{\natural^{}}_{\mathcal{H}}))$ 
with coefficients in $\mathsf{W}.$ 
\end{theorem}

We give an application of the above spectral sequence. Let $k$ be a field. Recall that a Hopf 
algebra $\mathcal{H}$ over $k$ is called semisimple
if it is semisimple as an algebra. It is shown in~\cite{sw69} that $\mathcal{H}$ is semisimple if and only if there is a right integral
$t \in \mathcal{H}$ with $\epsilon(t)=1$. Recall that a right integral in $\mathcal{H}$ is an element $t \in \mathcal{H}$ such
that for all $h \in \mathcal{H}, th=\epsilon(h) t.$ Now, it is easy to see that if $\mathcal{H}$ is semisimple, for any left
$\mathcal{H}$-module $M$, we have $H_0(\mathcal{H},M)=M_{\mathcal{H}}$ and $H_i(\mathcal{H},M)=0$ for $i>0$. In fact, we have the
following homotopy operator 
$h:  \mathcal{H}^{\otimes n}  \otimes M \rightarrow   \mathcal{H}^{\otimes (n+1)} \otimes M, n \ge 0,$
$$ h( h_1 \otimes \dots \otimes h_n \otimes m ) =  t \otimes h_1 \otimes \dots \otimes h_n  \otimes m.$$
One can check that $\delta h + h \delta =  id$.
Note that semisimple Hopf algebras have bijective antipodes.

\begin{prop}
Let $\mathcal{H}$ be a semisimple Hopf algebra. Then there is a natural isomorphism of cyclic homology groups
$$HC_{\bullet} (A \rtimes \mathcal{H}; \mathsf{W})= HC_{\bullet} (C^{\mathcal{H}}_{\bullet}(A); \mathsf{W}), $$
where $C^{\mathcal{H}}_{\bullet}(A)= H_0(\mathcal{H}, \mathsf{C}_{\bullet}(A^{\mathcal{H}}_{\natural}))$ is the cyclic module of equivariant chains.
\end{prop}
\begin{proof}
Since $\mathcal{H}$ is semisimple, we have $\mathsf{E}^1_{pq}=0$ for $q>0$ and the spectral sequence collapses. 
The first row of $\mathsf{E}^1$ is exactly 
$H_{0}(\mathcal{H},\mathsf{C}_{\bullet}(A^{\mathcal{H}}_{\natural}))=C_{\bullet}^{\mathcal{H}}(A).$
\end{proof}


\section{The Cocylindrical Module $\mathbf {C} \mathbf{\natural^{\textbf{}}} \mathbf{\mathcal{H}}.$}

Let $ \mathcal{H}$ be a Hopf algebra with a bijective antipode and $C$ a left $\mathcal{H}$-module coalgebra.
Our aim in this section is to establish a spectral sequence for computing the cyclic cohomology of the crossed
product coalgebra $C \! > \! \blacktriangleleft \! \mathcal{H}.$ To this end, we introduce
the cocylindrical module 
$$ C \natural^{\textbf{}} \mathcal{H} = \{ \mathcal{H}^{\otimes(p+1)}\otimes C^{\otimes (q+1)} \}_{p,q \ge 0}. $$  
We define the vertical and horizontal
operators $\tau_{p,q} ,\partial_{p,q} ,\sigma_{p,q}$ and $\bar{\tau}_{p,q} ,\bar{\partial}_{p,q} ,\bar{\sigma}_{p,q}$ 
as follows:
\begin{eqnarray} \label{eq:opcy5} \notag
& &\tau_{p,q}( g_0 , \dots , g_p \mid a_0 , \dots , a_q)   \\ \notag
& &= (g_0^{(1)}, \dots , g_p^{(1)} \mid  a_1, \dots ,a_{q}, (g_0^{(0)} S(g_0^{(2)})  
\dots g_p^{(0)} S(g_p^{(2)}))\cdot a_0),  \\ \notag
& &\sigma_{p,q}^i(g_0,\dots,g_p \mid a_0 , \dots , a_q)  \\\notag
& &= (g_0,\dots,g_p \mid a_0,\dots , a_{i}, a_{i+2},\dots , a_q) \epsilon(a_i), \;\;\; 0 \le i <q,  \\ 
& &\partial_{p,q}^i(g_0,\dots,g_p \mid a_0 , \dots , a_q) \\\notag
& &= (g_0,\dots,g_p \mid a_0,\dots , a_i^{(0)} , a_{i}^{(1)},\dots , a_q), \;\;\; 0 \le i \le q,\\\notag
& &\partial_{p,q}^{q+1}( g_0 , \dots , g_p \mid a_0 , \dots , a_q)   \\\notag
& &= (g_0^{(1)}, \dots , g_p^{(1)} \mid a_0^{(1)}, a_1, \dots ,a_{q},(g_0^{(0)} S(g_0^{(2)}) 
\dots g_p^{(0)} S(g_p^{(2)})) \cdot a_0^{(0)}),  \notag
\end{eqnarray}
where $g \cdot a$ denotes the action of $g \in \mathcal{H}$ on $a \in C,$
\begin{eqnarray} \label{eq:opcy6} \notag
& &\bar{\tau}_{p,q}( g_0 , \dots , g_p \mid a_0 , \dots , a_q)
=(g_1, \dots , g_{p},g_0^{(1)} \mid S^{-1}(g_0^{(0)}S(g_0^{(2)})) \cdot (a_0, \dots a_{q}) ),\\ \notag
& &\bar{\sigma}_{p,q}^i(g_0,\dots,g_p \mid a_0 , \dots , a_q)\\ \notag
&&\hspace{2cm}=(g_0,\dots,g_i, g_{i+2},\dots,g_p \mid a_0,\dots , a_q) \epsilon(g_i), \;\;\; 0 \le i < p,  \\ 
& &\bar{\partial}_{p,q}^i (g_0,\dots,g_p \mid a_0 , \dots , a_q)\\ \notag
&&\hspace{2cm}=(g_0,\dots,g_i^{(0)},g_{i}^{(1)},\dots,g_p \mid a_0,\dots , a_q), \;\;\; 0 \le i \le p, \\ \notag
& &\bar{\partial}_{p,q}^{p+1}( g_0 , \dots , g_p \mid a_0 , \dots , a_q)
=(g_0^{(3)},g_1, \dots , g_{p},g_0^{(1)} \mid S^{-1}(g_0^{(0)}S(g_0^{(2)})) \cdot (a_0, \dots a_{q})), \notag
\end{eqnarray} 
where $g \cdot (a_0, \dots,a_q)=(g^{(0)} \cdot a_0,\dots,g^{(q)} \cdot a_q)$ denotes the diagonal action of $\mathcal{H}$
on $C^{\otimes n}$.

The proof of the following theorem is similar to the proof of Theorem 2.1 and we omit it.
\begin{theorem}
Endowed with the operators defined in (\ref{eq:opcy5}),(\ref{eq:opcy6}), $C \natural^{\text{}} \mathcal{H}$ is a 
cocylindrical module.  
\end{theorem}

\begin{corollary}$\Delta (C \natural^{\textbf{}} \mathcal{H})$ is a cocyclic module.
\end{corollary}

\begin{theorem} \label{iso2}
 We have an isomorphism of cocyclic modules 
\[\Delta (C \natural^{} \mathcal{H}) \cong \mathsf{C}^{\bullet}(C  \!> \! \blacktriangleleft \!\mathcal{H}).\]
\end{theorem}
\begin{proof}
One can check that $\phi$ and $\psi$ define isomorphisms, inverse to each other, between
$\Delta (C \natural \mathcal{H})$ and $\mathsf{C}^{\bullet}(C  \!> \! \blacktriangleleft \!\mathcal{H})$ where,
\begin{eqnarray*}
& &\phi ( a_0 \otimes g_0,\dots ,a_n \otimes g_n )=
(g_0^{(n)} , 
g_1^{(n-1)}, \dots ,g_{n-1}^{(1)},g_n \mid \\
& & \hspace{0.5 cm}a_0,(g_0^{(0)}S(g_0^{(2n)})) \cdot a_1 , \dots ,
(g_{n-1}^{(0)}S(g_{n-1}^{(2)}) \dots\\
& & \hspace{4cm} g_1^{(n-2)}S(g_{1}^{(n)})g_0^{(n-1)}S(g_0^{(n+1)})) \cdot a_n),  
\end{eqnarray*}
and
\begin{eqnarray*}
& &\psi (g_0,g_1, \dots ,g_{n-1}, g_n \mid a_0 , \dots ,a_n)=
(a_0 \otimes g_0^{(n)}, S^{-1}(g_0^{(n-1)}S(g_0^{(n+1)})) \cdot a_1  \otimes  
 g_1^{(n-1)},\\
& &S^{-1}(g_1^{(n-2)}S(g_1^{(n+2)})g_0^{(n-2)}S(g_0^{(n+2)})) \cdot a_1  \otimes  
 g_1^{(n-1)},
\dots,\\
& &\hspace{2cm}  S^{-1}(g_{n-1}^{(0)}S(g_{n-1}^{(2)}) 
\dots g_1^{(0)}S(g_{1}^{(2n-2)})g_0^{(0)}S(g_0^{(2n)})) \cdot a_n \otimes g_n). 
\end{eqnarray*}
\end{proof}

We define a left $\mathcal{H}$-coaction  on the first column of 
$C \natural^{} \mathcal{H}$, $C^{\natural^ {}}_{\mathcal{H}}$ 
= $ \{ \mathcal{H} \otimes C^{\otimes (n+1)} \}_{n \ge 0}$, by 
\begin{equation} \label{eq:ac}
\boldsymbol{\Delta} (g \mid a_0, \dots , a_n ) = ( S(g^{(4)})g^{(1)} \mid g^{(3)} \mid (g^{(2)}S^{-1}(g^{(0)})) \cdot 
( a_0, \dots ,a_n)). 
\end{equation}

We define $\mathsf{C}^q(C_{\mathcal{H}}^{\natural^{}})=\mathcal{H} \otimes C^{\otimes (q+1)}$ and
$\mathsf{C}^p(\mathcal{H},\mathsf{C}^q(C_{\mathcal{H}}^{\natural^{}}))= \mathcal{H}^{\otimes p} \otimes 
\mathsf{C}^q(C_{\mathcal{H}}^{\natural^{}}).$ So we can construct $H^p(\mathcal{H},\mathsf{C}^q(C^{\natural^{}}_{\mathcal{H}}))$,
the 
cohomology of the coalgebra $\mathcal{H}$ with coefficients in the 
comodule $\mathsf{C}^q(C_{\mathcal{H}}^{\natural^{}}).$       
We define $C^{n}_{\mathcal{H}^{}}(C)$ as the coinvariant space of $\mathcal{H} \otimes C^{\otimes (n+1)}$ under 
the above
coaction, so that $C^{n}_{\mathcal{H}^{}}(C)$ is the space of all $(g \mid a_0, \dots , a_n ) $ such that \[
\boldsymbol{\Delta} (g \mid a_0, \dots , a_n ) = (1 \mid g \mid a_0, \dots , a_n ). \]

Now we replace the cocylindrical module $(C \natural^{} 
\mathcal{H},(\partial,\sigma,\tau),(\bar{\partial},\bar{\sigma},\bar{\tau}))$ with an isomorphic cocylindrical module
$(\mathsf{C}(\mathcal{H},\mathsf{C} ( C^{\natural^{}}_{\mathcal{H}})),( \mathfrak{d},\mathfrak{s},\mathfrak{t}),(\bar{\mathfrak{d}},\bar{\mathfrak{s}},\bar{\mathfrak{t}}))$,
under the isomorphism defined by 
the maps $\boldsymbol{\beta} : (C \natural^{} \mathcal{H})_{p,q} \rightarrow \mathsf{C}^p(\mathcal{H},\mathsf{C}^q(C_{\mathcal{H}}^{\natural^{}}))$
and $\boldsymbol{\gamma} :  
\mathsf{C}^p(\mathcal{H},\mathsf{C}^q(C_{\mathcal{H}}^{\natural^{}})) \rightarrow  (C \natural^{} \mathcal{H})_{p,q},$
\begin{eqnarray*}
& &\boldsymbol{\beta} ( g_0, \dots , g_p \mid a_0 , \dots , a_q ) = 
( S(g_{0}^{(1)}) \cdot (g_1, \dots , g_p )\mid g_0^{(0)} \mid
a_0, \dots , a_q ), \\
& &\boldsymbol{\gamma} ( g_1 , \dots , g_p \mid g \mid a_0 , \dots , a_q) = 
(g^{(0)},g^{(1)} \cdot (g_1 , \dots , g_p)  \mid 
a_0, \dots , a_q ).
\end{eqnarray*}

One can check that $\boldsymbol{\beta} \circ \boldsymbol{\gamma} = \boldsymbol{\gamma} \circ \boldsymbol{\beta} = id $. 
We find the operators
$( \mathfrak{d},\mathfrak{s},\mathfrak{t}),(\bar{\mathfrak{d}},\bar{\mathfrak{s}},\bar{\mathfrak{t}})$
under this transformation. Since the computations are similar to the comodule algebra case, we only give the final results. \\

$ \hspace{0.2cm} \bullet \quad $ 
For $\bar{\mathfrak{b}} = \boldsymbol{\beta} \bar{b} \boldsymbol{\gamma}$, we have 
\begin{eqnarray*}
& &\bar{\mathfrak{b}}( g_1 , \dots , g_p \mid g \mid a_0 , \dots , a_q) =
\boldsymbol{\delta}(g_1, \dots ,g_p) \mid g \mid a_0 , \dots , a_q ). 
\end{eqnarray*}
Here $\boldsymbol{\delta}$ is the Hopf-comodule coboundary defined by
\begin{eqnarray} 
& &\boldsymbol{\delta} (g_1,g_2,\dots , g_p,m)= (1,g_1,\dots ,g_p,m) +\\
& &\sum_{i=1}^{p} (-1)^i (g_1 , \dots ,g_i^{(0)}, g_{i}^{(1)},\dots , g_p ,m) + (-1)^{p+1} 
(g_1,\dots,g_{p},\Delta_{\tiny{M}} (m)), \notag
\end{eqnarray}
where $m \in M$ and $M$ is a left $\mathcal{H}$-comodule with structure map $\Delta_{\tiny{M}}$.
So we have,
\begin{eqnarray*}
& &\bar{\mathfrak{d}}^0 (g_1, \dots ,g_p \mid g \mid a_0, \dots , a_q) = (1,g_1 , 
\dots , g_p \mid g \mid a_0 , \dots , a_q),\\
& &\bar{\mathfrak{d}}^i (g_1, \dots ,g_p \mid g \mid a_0, \dots , a_q) = 
(g_1 , \dots , g_i^{(0)}, g_i^{(1)} ,\dots , g_p \mid g \mid a_0 , \dots , a_q),\;\; 1 \le i \le p,\\
& &\bar{\mathfrak{d}}^{p+1} (g_1, \dots ,g_p \mid g \mid a_0, \dots , a_q) = (g_1 , \dots , 
g_{p} \mid \boldsymbol{\Delta}( g \mid a_0 , \dots ,  a_q)).
\end{eqnarray*}

$ \hspace{0.2cm} \bullet \quad $ For $\bar{\mathfrak{t}}= 
\boldsymbol{\beta} \bar{\tau} \boldsymbol{\gamma}$ we have,
\begin{eqnarray*}
& &\bar{\mathfrak{t}} (g_1 , \dots , g_p \mid g \mid a_0, \dots , a_q) =\\
&=&( S(g_1^{(2)})\cdot (g_2 \dots ,g_p),S(g^{(4)} g_1^{(1)}) g^{(1)}\mid g^{(3)} g_1^{(0)} \mid 
S^{-1}(g^{(0)}S(g^{(2)})) \cdot (a_0, \dots ,a_q)).
\end{eqnarray*}
For  $\bar{\mathfrak{s}}^i= \boldsymbol{\beta} \bar{\sigma}^i \boldsymbol{\gamma},$ $\;\; 0 \le i < p, \;$ we have
\begin{eqnarray*}
& &\bar{\mathfrak{s}}^i (g_1 , \dots , g_p \mid g \mid a_0, \dots , a_q) 
= (g_1 , \dots ,g_i,g_{i+2}, \dots, g_p \mid g \mid a_0, \dots , a_q) \epsilon(g_{i+1}).
\end{eqnarray*}
$ \hspace{0.2cm} \bullet \quad $ Next we compute the operator  $\mathfrak{b} = \boldsymbol{\beta} b \boldsymbol{\gamma}$. We have 
\begin{eqnarray*}
& &\mathfrak{b}(g_1, \dots , g_p \mid g \mid a_0 , \dots , a_q) = \boldsymbol{\beta} b (g^{(0)},g^{(1)} \cdot
(g_1, \dots ,g_p) \mid a_0 , \dots , a_q )\\
&=& (\sum_{0 \le i \le q} (-1)^i (g_1, \dots ,g_p \mid g  \mid a_0 , \dots ,a_i^{(0)},a_i^{(1)},\dots, a_q )\\
&+& (-1)^{p+1}
( g_1^{(1)}, \dots ,  g_p^{(1)} \mid g^{(1)} \mid a_0^{(1)}, a_1, \dots ,a_q, \\
& & \hspace{6cm} (g^{(0)} g_1^{(0)} S(g_1^{(2)} ) \dots g_p^{(0)} S(g_p^{(2)}) S(g^{(2)})) \cdot a_0^{(0)}). 
\end{eqnarray*}
So we conclude that,
\begin{eqnarray*}
& &\mathfrak{d}^i(g_1,\dots,g_p \mid g \mid a_0,\dots,a_q)=(g_1 , 
\dots , g_p \mid g \mid a_0, \dots,a_{i}^{(0)}, a_{i}^{(1)},\dots,a_{q}),\;\; 0 \le i \le q, \\
& &\mathfrak{d}^{q+1}(g_1,\dots,g_p \mid g \mid a_0,\dots,a_q)= ( g_1^{(1)}, \dots ,  g_p^{(1)} \mid g^{(1)} \mid a_0^{(1)}, a_1, \dots ,a_q, \\
& &\hspace{6cm} (g^{(0)} g_1^{(0)} S(g_1^{(2)} ) \dots g_p^{(0)} S(g_p^{(2)}) S(g^{(2)})) \cdot a_0^{(0)}). 
\end{eqnarray*}
$ \hspace{0.2cm} \bullet \quad $ We consider  $\mathfrak{s}^i = \boldsymbol{\beta} \sigma^i 
\boldsymbol{\gamma}$ and $\mathfrak{t} = \boldsymbol{\beta} \tau \boldsymbol{\gamma}$. We have 
\begin{eqnarray*}
& &\mathfrak{s}^i (g_0,\dots,g_p \mid g \mid a_0, \dots ,a_q)
=(g_1, \dots ,g_p \mid g \mid a_0 , \dots ,a_i,a_{i+2},\dots, a_q ) \epsilon{(a_{i+1})}. 
\end{eqnarray*}
Also,
\begin{eqnarray*}
& &\mathfrak{t}(g_1,\dots,g_p \mid g \mid a_0,\dots,a_q)= ( g_1^{(1)}, \dots ,  g_p^{(1)} \mid g^{(1)} \mid a_1, \dots ,a_q, \\
& &\hspace{6cm} (g^{(0)} g_1^{(0)} S(g_1^{(2)} ) \dots g_p^{(0)} S(g_p^{(2)}) S(g^{(2)})) \cdot a_0 ).  
\end{eqnarray*}

By the above computations we can state the following theorems.
\begin{theorem} \label{ta:iso1}
The complex $(\mathsf{C} (C \natural^{} \mathcal{H}) \boxtimes \mathsf{W}, b + 
\mathbf{u} B , \bar{b} + \mathbf{u} \bar{B}) $ 
is isomorpaic to the complex 
$(\mathsf{C}(\mathcal{H}, \mathsf{C}(C^{\natural^{}}_{\mathcal{H}} \boxtimes \mathsf{W})), 
\mathfrak{b} + \mathbf{u} \mathfrak{B} , \bar{\mathfrak{b}} + \mathbf{u} \bar{\mathfrak{B}}), $
where $\mathfrak{\bar{b}}$ is the Hopf-comodule coboundary. 
\end{theorem}
\begin{theorem} \label{ta:iso2}
For  $p \ge 0$, $H^p(\mathcal{H},\mathsf{C}^{\bullet} (C^{\natural^{}}_{\mathcal{H}}))$, with the operators defined as follows,
are cocyclic modules: 
\begin {eqnarray*} \label{eq:equ}
& &\mathfrak{t}(g_1,\dots,g_p \mid g \mid a_0,\dots,a_q)= ( g_1^{(1)}, \dots ,  g_p^{(1)} 
\mid g^{(1)} \mid a_1, \dots ,a_q, \\
& &\hspace{6cm} (g^{(0)} g_1^{(0)} S(g_1^{(2)} ) \dots g_p^{(0)} S(g_p^{(2)}) S(g^{(2)})) \cdot a_0 ),\\  
& &\mathfrak{d}^i(g_1,\dots,g_p \mid g \mid a_0,\dots,a_q)=(g_1 , \dots , g_p \mid g \mid a_0, 
\dots,a_{i}^{(0)}, a_{i}^{(1)},\dots,a_{q}),\;\; 0 \le i \le q, \\
& &\mathfrak{d}^{q+1}(g_1,\dots,g_p \mid g \mid a_0,\dots,a_q)= ( g_1^{(1)}, \dots ,  
g_p^{(1)} \mid g^{(1)} \mid a_0^{(1)}, a_1, \dots ,a_q, \\
& &\hspace{6cm} (g^{(0)} g_1^{(0)} S(g_1^{(2)} ) \dots g_p^{(0)} S(g_p^{(2)}) S(g^{(2)})) \cdot a_0^{(0)}),\\  
& &\mathfrak{s}^i(g_1,\dots,g_p \mid g \mid a_0,\dots,a_q)=(g_1, \dots , g_p \mid g \mid  a_0 , 
\dots ,a_i, a_{i+2},\dots, a_q) \epsilon(a_{i+1}).
\end{eqnarray*}
\end{theorem}

\begin{proof}
We prove this for $p=0,$ the general case being similar. We have 
$H^0(\mathcal{H},\mathsf{C}^{\bullet} (C^{\natural^{}}_{\mathcal{H}})) \cong C^{\bullet}_{\mathcal{H}}(C).$
The operators are well defined on the coinvariant space. For example for 
$\mathfrak{t},$ since $ \mathfrak{t} \bar{\mathfrak{b}} = \bar{\mathfrak{b}} \mathfrak{t} $
and  $\bar{\mathfrak{b}} (g \mid a_0, \dots , a_n ) = \boldsymbol{\Delta} (g \mid a_0, \dots , a_n ) - 
(1 \mid g \mid a_0, \dots , a_n ),$ if $ \mathfrak{t} \bar{\mathfrak{b}} = 0$
then $\bar{\mathfrak{b}} \mathfrak{t} = 0$ .

To show that $ \mathfrak{t} ^{n+1} = 1,$ we see that in the coinvariant space we have,
\begin{eqnarray*}
(1 \mid g \mid a_0, \dots , a_n ) \equiv ( S(g^{(4)})g^{(1)} \mid g^{(3)} \mid (g^{(2)}S^{-1}(g^{(0)})) 
\cdot ( a_0, \dots ,a_n))  
\end{eqnarray*}
So we have, 
\begin{eqnarray*}
\bar{\mathfrak{s}}^0 \circ \bar{\mathfrak{t}}(1 \mid g \mid a_0, \dots , a_n ) \equiv ( S(g^{(4)})g^{(1)} \mid g^{(3)} \mid (g^{(2)}S^{-1}(g^{(0)})) 
\cdot ( a_0, \dots ,a_q)), 
\end{eqnarray*}
where,
$$\bar{\mathfrak{s}}^0 \circ \bar{\mathfrak{t}} (1 \mid g \mid a_0, \dots , a_n ) =(g \mid a_0, \dots , a_n ),$$
and,
\begin{eqnarray*}
& &\bar{\mathfrak{s}}^0 \circ \bar{\mathfrak{t}}( S(g^{(4)})g^{(1)} \mid g^{(3)} \mid (g^{(2)}S^{-1}(g^{(0)})) 
\cdot ( a_0, \dots ,a_q)) =\\
& &\bar{\mathfrak{s}}^0 (S(g^{(8)} S(g^{(9)}) g^{(2)}) g^{(5)} \mid g^{(7)} S(g^{(10)}) g^{(1)} 
\mid \\
& & \hspace{2cm}S^{-1}(g^{(4)}S(g^{(6)}))S^{-1}(g^{(0)}S(g^{(3)}))\cdot( a_0, \dots , a_n ))\\
&=&( g^{(1)} \mid S^{-1}(g^{(0)}S(g^{(2)}))\cdot( a_0, \dots , a_n )).
\end{eqnarray*}
Therefore in the coinvariant space we have \[(g \mid a_0, \dots , a_n ) 
\equiv( g^{(1)} \mid S^{-1}(g^{(0)} S(g^{({2})})) \cdot( a_0, \dots , a_n )), \]
and now it is easy to check that $\mathfrak{t}^{n+1}=1$.
\end{proof}

We use again the Eilenberg-Zilber theorem for cocylindrical modules combined with Theorem~\ref{iso2}
 to conclude 
$$H^{\bullet}(Tot(C \natural^{} \mathcal{H}) \boxtimes \mathsf{W}) \cong HC^{\bullet} (\Delta( C \natural^{}  \mathcal{H});
\mathsf{W}) \cong HC^{\bullet} (C \! > \! \blacktriangleleft \! \mathcal{H};
\mathsf{W}).$$
To compute the cohomology of the mixed complex \\
$(Tot(C \natural^{} \mathcal{H}), b + \bar{b} + \mathbf{u} (B + \bar{B})),$
we filter it by the subspaces \[
\mathsf{F}^i Tot^{n}((C \natural \mathcal{H} ) \boxtimes \mathsf{W}) = \sum_{q \ge i,\;p+q=n} (\mathcal{H}^{\otimes(p+1)} \otimes C^{\otimes(q+1)}) \boxtimes \mathsf{W}. \]
If again we separate the operator 
$b+ \bar{b} + \mathbf{u}(B + T \bar{B})$ as $\bar{b}+(b+  \mathbf{u}B) + \mathbf{u}T \bar{B}$, 
from  Theorems (\ref{th:iso1}),(\ref{th:iso2}) we can deduce the 
following theorem. 

\begin{theorem}
There is an spectral sequence that converges to $HC^{\bullet} (C \! > \! \blacktriangleleft \! \mathcal{H}; \mathsf{W})$. The $\mathsf{E}_0$-term of this spectral sequence is isomorphic to the complex 
\[ \mathsf{E}_0^{pq}=(\mathsf{C}^p(\mathcal{H},\mathsf{C}^q(C ^{\natural^{}}_{\mathcal{H}}) \boxtimes \mathsf{W} ),
\boldsymbol{\delta)}, \]
and the $\mathsf{E}_1$-term is  \[
\mathsf{E}_1^{pq}= (H^p( \mathcal{H} , \mathsf{C}^q(C ^{\natural^{}}_{\mathcal{H}}) \boxtimes \mathsf{W} ) , 
\mathfrak{b} + \mathbf{u} \mathfrak{B}). \]
The $\mathsf{E}_2$-term of the spectral sequence is
\[ \mathsf{E}_2^{pq} = HC^p( H^q(\mathcal{H}, \mathsf{C}^{\bullet}(C ^{\natural^{}}_{\mathcal{H}}));\mathsf{W}), \]
the cyclic cohomology of the cyclic comodules $H^p(\mathcal{H},\mathsf{C}^{\bullet}(C ^{\natural^{}}_{\mathcal{H}}))$ with coefficients in $\mathsf{W}.$ 
\end{theorem}

We give an application of the above spectral sequence. Let $k$ be a field. Recall that a Hopf algebra $\mathcal{H}$ over
$k$ is called cosemisimple if $\mathcal{H}$ is cosemisimple as a coalgebra, that is, every left $\mathcal{H}$-comodule is
completely reducible~\cite{sw69}. One knows that a Hopf algebra is cosemisimple if and only if there exists a left integral
$x \in \mathcal{H^{\ast}}$ with $x(1)=1$~(\cite{sw69}, Theorem 14.0.3). It is easy to see that if $\mathcal{H}$ is 
cosemisimple and $M$ is an $\mathcal{H}$-bicomodule, then the coalgebra (Hochschild) cohomology groups satisfy
$H^i(\mathcal{H},M)=0$ for $i > 0,$ and $H^0(\mathcal{H},M)= M^{co\mathcal{H}},$ the subspace of coinvariants of 
the bicomodule $M$. In fact, we have the following homotopy 
operator $h: \mathcal{H}^{\otimes n} \otimes M \rightarrow \mathcal{H}^{\otimes(n-1)} \otimes M, n \ge 1,$
\[
h(g_1,\dots,g_n,m)=x(g_1)(g_2,\dots,g_n,m).
\]
One can check theat $\delta h + h \delta=id.$ Note that the antipode of $\mathcal{H}$ is bijective if $\mathcal{H}$ is 
cosemisimple.
\begin{prop}
Let $\mathcal{H}$ be a cosemisimple Hopf algebra. Then there is a natural isomorphism of cyclic and Hochschild cohomology
groups 
\[ HC^{\bullet}(C  \!> \! \blacktriangleleft \!\mathcal{H}) \simeq HC^{\bullet}(\mathsf{C}^{\bullet}_{\mathcal{H}}(C)),\]
\[ HH^{\bullet}(C  \!> \! \blacktriangleleft \!\mathcal{H}) \simeq HH^{\bullet}(\mathsf{C}^{\bullet}_{\mathcal{H}}(C)).\]
\end{prop}
\begin{proof}
Since $\mathcal{H}$ is cosemisimple, we have $\mathsf{E}^{p,q}_1=0$ for $p > 0$ and the spectral sequence collapses. 
The first column
of $\mathsf{E}_1$ is exactly $H^{\bullet}(\mathcal{H}, C^{\natural}_{\mathcal{H}})= \mathsf{C}^{\bullet}_{\mathcal{H}}(C)$.
\end{proof}

\section*{Acknowledgments}

We would like to thank the referee for pointing out a number of typographical errors in the original version of this paper.


\end{document}